\documentclass[12pt]{article}
\usepackage{mathrsfs}
\usepackage{epic,eepic,epsf,epsfig}
\usepackage{amsfonts,srcltx,mathrsfs}
\usepackage{multirow}
\usepackage{amsmath}
\usepackage{amssymb}
\usepackage{amsbsy}
\usepackage{graphicx}
\usepackage{amsfonts}%,srcltx}
\usepackage{color}
\usepackage{setspace}

\newtheorem{lem}{Lemma}[section]%
\newtheorem{theorem}[lem]{Theorem}%
\newtheorem{cor}[lem]{Corollary}%

\def\nd{\mathrel{\bigm|\kern-.7em/}}

\def\f{\noindent}

\def\P\GammaL{\hbox{\rm P\GammaL}}

\def\mod{\hbox{\rm mod }}

\begin{document}
\title{Spectral extrema of graphs forbidding a fan}

\footnotetext{E-mails: zhangwq@pku.edu.cn}

\author{Wenqian Zhang\\
{\small School of Mathematics and Statistics, Shandong University of Technology}\\
{\small Zibo, Shandong 255000, P.R. China}}
\date{}
\maketitle

\begin{abstract}
For a graph $G$, its spectral radius is the largest eigenvalue of its adjacency matrix. A fan $H_{\ell}$ is a graph obtained by connecting a single vertex to all vertices of a path of order $\ell\geq4$. Let ${\rm SPEX(n,H_{\ell})}$ be the set of all extremal graphs $G$ of order $n$ with the maximum spectral radius, where $G$ contains no $H_{\ell}$ as a subgraph. In this paper, we completely characterized the graphs in ${\rm SPEX(n,H_{\ell})}$ for any $\ell\geq4$ and sufficiently large $n$. An interesting phenomenon was revealed: ${\rm SPEX(n,H_{2k+2})}\subseteq {\rm SPEX(n,H_{2k+3})}$ for any $k\geq1$ and sufficiently large $n$.

\bigskip

\f {\bf Keywords:} spectral radius; spectral extremal graph; fan; path.

\end{abstract}

 \baselineskip 17 pt

\section{Introduction}

All graphs considered in this paper are finite and undirected.
For a graph $G$, let $\overline{G}$ be its complement. The vertex set and edge set of $G$ are denoted by $V(G)$ and $E(G)$, respectively. Let $e(G)=|E(G)|$. For a vertex $u$, let $d_{G}(u)$ be its degree. Let $\delta(G)$ or $\Delta(G)$ denote the minimum or maximum degree of $G$. For any $S\subseteq V(G)$, let $G[S]$ be the  subgraph of $G$ induced by $S$, and let $G-S=G[V(G)-S]$. For two vertices $u$ and $v$, we say that $u$ is a neighbor of $v$ or $u\sim v$, if they are adjacent in $G$. Let $N_{S}(u)$ be the set of neighbors in $S$ of $u$, and let $d_{S}(u)=|N_{S}(u)|$. For two disjoint subsets $S,T\subseteq V(G)$, let $e_{G}(S,T)$ be the number of edges between $S$ and $T$ in $G$.     For two graphs $G_{1}$ and $G_{2}$, let $G_{1}\vee G_{2}$ be their join, which is obtained from their disjoint union $G_{1}\cup G_{2}$, by connecting each vertex in $G_{1}$ to all the vertices in $G_{2}$.  For a certain integer $n$, let $K_{n},P_{n}$ and $K_{1,n-1}$ be the complete graph (or a clique), the path and the star graph of order $n$, respectively. For $r\geq2$, let $T(n,r)$ be the Tur\'{a}n graph of order $n$ with $r$ parts. For any terminology used but not defined here, one may refer to \cite{CRS}.

 Let $G$ be a graph with vertices $v_{1},v_{2},...,v_{n}$. The {\em adjacency matrix} of $G$ is  $A(G)=(a_{ij})_{n\times n}$, where $a_{ij}=1$ if  $v_{i}\sim v_{j}$, and $a_{ij}=0$ otherwise. The {\em spectral radius} $\rho(G)$ of $G$ is the largest eigenvalue of $A(G)$. By the Perron--Frobenius theorem,  $\rho(G)$ has a non-negative eigenvector (called Perron vector), and has a positive eigenvector if $G$ is connected.
For a set of graphs $\mathcal{F}$, a graph $G$ is call $\mathcal{F}$-free if $G$ does not contain any member in $\mathcal{F}$  as a subgraph. Let ${\rm SPEX}(n,F)$ be the set of $\mathcal{F}$-free graphs of order $n$ with the maximum spectral radius. We also use $F$ instead of $\mathcal{F}$ when $\mathcal{F}=\left\{F\right\}$.

 In 2010, Nikiforov \cite{N1} proposed a spectral version of the Tur\'{a}n-type
problem: what is the maximum spectral radius of an $F$-free graph of order $n$?
 In recent years, this problem has been studied for many kinds of $F$ (see \cite{BDT,CFTZ,CDT1,CDT2,CDT3,DKLNTW,FZ,FTZ,FL,LL,LP1,LP2,LZZ,LM,MLZ,
 NWK,N,N1,N2,WKX,ZHL,ZL,ZL1,ZLX,Z1,Z2}).

 For $\ell\geq4$, the fan graph is defined as $H_{\ell}=K_{1}\vee P_{\ell}$. Recently, many researchers concerned the spectral extrema of $H_{\ell}$-free graphs. For $H_{\ell}$-free graphs $G$ with $e(G)$ fixed and $|G|$ released, the spectral extrema of $G$ is studied (see \cite{CY,GL,LiZZ,LZG,YLP,ZW}). 
 For an integer $k\geq1$, $G$ is called {\em nearly $k$-regular} if all its vertices have degree $k$ except one vertex with degree $k-1$. Clearly, if $G$ is nearly $k$-regular, then $k|G|$ is odd. (Note that $H_{4}$ is the square of $P_{5}$.) Zhao and Park \cite{ZP} characterized the graphs in ${\rm SPEX}(n,H_{4})$ 
 
 \begin{theorem}{\rm (\cite{ZP})}\label{even 1}
For $n\geq6$, the unique graph $G$ in $\rm{SPEX}(n,H_{4})$ is obtained from
 a complete bipartite graph with parts $L$ and $R$, by embedding a (nearly) $1$-regular graph in $G[L]$. Moreover,\\
\begin{equation}
|L|=\left\{
\begin{array}{lr}
\frac{n}{2},&n\equiv 0(\mod4);\\
\\
\frac{n-1}{2},&n\equiv 1(\mod4);\\
\\
\frac{n}{2},&n\equiv 2(\mod4);\\
\\
\frac{n+1}{2},&n\equiv 3(\mod4).\\
\end{array}
\right.\notag
\end{equation}.
\end{theorem}

For $k\geq3$, Yuan, Liu and Yuan \cite{YLY} characterized the graphs in ${\rm SPEX}(n,H_{2k})$. 

\begin{theorem}{\rm (\cite{YLY})}\label{even 2}
For $k\geq3$ and sufficiently large $n$, any graph $G$ in $\rm{SPEX}(n,H_{2k})$ is obtained from
 a complete bipartite graph with parts $L$ and $R$, by embedding a (nearly) $(k-1)$-regular $P_{2k}$-free graph in $G[L]$. Moreover, $\frac{n}{2}-1\leq |L|,|R|\leq\frac{n}{2}+1$.
\end{theorem}

 In this paper, we completely characterize the graphs in $\rm{SPEX}(n,H_{2k+3})$ for $k\geq1$ and large $n$.

\begin{theorem}\label{main1}
For $k\geq1$ and sufficiently large $n$, any graph $G$ in $\rm{SPEX}(n,H_{2k+3})$ is obtained from
 a complete bipartite graph with parts $L$ and $R$, by embedding a (nearly) $k$-regular $P_{2k+3}$-free graph in $G[L]$. Moreover, $|L|=\lfloor\frac{n}{2}\rfloor$ or $\lceil\frac{n}{2}\rceil$ for even $k\geq2$, and for odd $k\geq1$\\
\begin{equation}
|L|=\left\{
\begin{array}{lr}
\frac{n}{2},&n\equiv 0(\mod4);\\
\\
\frac{n-1}{2},&n\equiv 1(\mod4);\\
\\
\frac{n}{2},&k=1~and~n\equiv 2(\mod4);\\
\\
\frac{n}{2}-1~or~\frac{n}{2}+1,&k\geq3~and~n\equiv 2(\mod4);\\
\\
\frac{n+1}{2},&n\equiv 3(\mod4).\\
\end{array}
\right.\notag
\end{equation}.
\end{theorem}

Note that the spectral extremal graphs in Theorem \ref{main1} are completely determined. For example, when $k\geq2$ is even, the resulting graph $G$ has the same spectral radius whenever $|L|=\lfloor\frac{n}{2}\rfloor$ or $\lceil\frac{n}{2}\rceil$. When $k\geq3$ is odd and $n\equiv 2(\mod4)$, the resulting graph $G$ has the same spectral radius whenever $|L|=\frac{n}{2}-1$ or $\frac{n}{2}+1$. Furthermore, one can see that $G[L]$ is $k$-regular except the case of $k=1$ and $n\equiv 2(\mod4)$. By a similar proof as Theorem \ref{main1}, we can slightly refine Theorem \ref{even 2} as follows.

\begin{theorem}\label{main2}
For $k\geq1$ and sufficiently large $n$, any graph $G$ in $\rm{SPEX}(n,H_{2k+2})$ is obtained from
 a complete bipartite graph with parts $L$ and $R$, by embedding a (nearly) $k$-regular $P_{2k+2}$-free graph in $G[L]$. Moreover, $|L|=\lfloor\frac{n}{2}\rfloor$ or $\lceil\frac{n}{2}\rceil$ for even $k\geq2$, and for odd $k\geq1$\\
\begin{equation}
|L|=\left\{
\begin{array}{lr}
\frac{n}{2},&n\equiv 0(\mod4);\\
\\
\frac{n-1}{2},&n\equiv 1(\mod4);\\
\\
\frac{n}{2},&k=1~and~n\equiv 2(\mod4);\\
\\
\frac{n}{2}-1~or~\frac{n}{2}+1,&k\geq3~and~n\equiv 2(\mod4);\\
\\
\frac{n+1}{2},&n\equiv 3(\mod4).\\
\end{array}
\right.\notag
\end{equation}.
\end{theorem}

From Theorem \ref{main1} and Theorem \ref{main2}, we see $\rm{SPEX}(n,H_{2k+2})\subseteq \rm{SPEX}(n,H_{2k+3})$ for any $k\geq1$ and large $n$. Thus, Theorem \ref{main1} essentially strengthens Theorem \ref{main2}.

The rest of this paper is organized as follows. In Section 2, we include some lemmas, which will be used in the proof of Theorem \ref{main1}.  In Section 3, we give a general result for spectral extremal graphs. In Section 4, we study spectral radius using walks in graphs. In Section 5, we give the proof of Theorem \ref{main1}.

\section{Preliminaries}

 To prove Theorem \ref{main1}, we first include several lemmas.

\begin{lem}{\rm (\cite{CRS})}\label{subgraph}
If $H$ is a subgraph of a connected graph $G$, then $\rho(H)\leq\rho(G)$, with equality if and only if $H=G$.
\end{lem}

The following lemma is a variation (with a very similar proof) of Theorem 8.1.3 in \cite{CRS}.

\begin{lem}{\rm (\cite{CRS})}\label{eigenvector trans}
Let $G$ be a connected graph with a Perron vector $\mathbf{x}=(x_{w})_{w\in V(G)}$. For a vertex $u\in V(G)$, Let $G'$ be the graph obtained from $G$ by deleting edges $uv_{1},uv_{2},...,uv_{s}$, and adding edges $uw_{1},uw_{2},...,uw_{t}$, where $s,t\geq1$. If $\sum_{1\leq j\leq t}x_{w_{j}}\geq\sum_{1\leq i\leq s}x_{v_{i}}$ and $\left\{v_{1},v_{2},...,v_{s}\right\}\neq\left\{w_{1},w_{2},...,w_{s}\right\}$, then $\rho(G')>\rho(G)$.
\end{lem}

\medskip

The following is the Spectral Stability Lemma due to Nikiforov \cite{N2}.

\begin{lem} {\rm (\cite{N2})}\label{spec-stability}
Suppose $r\geq2,\frac{1}{\ln n}<c<r^{-8(r+1)(r+21)}$ and $0<\epsilon<2^{-36}r^{-24}$. Let $G$ be a graph of order $n$. If $\rho(G)>(\frac{r-1}{r}-\epsilon)n$, then one of the following holds:\\
$(i)$ $G$ contains a complete $(r+1)$-partite graph $K_{\lfloor c\ln n\rfloor,\lfloor c\ln n\rfloor,...,\lfloor c\ln n\rfloor,\lceil n^{1-\sqrt{c}}\rceil}$;\\
$(ii)$ $G$  differs from $T(n,r)$ in fewer than $(\epsilon^{\frac{1}{4}}+c^{\frac{1}{8r+8}})n^{2}$ edges.
\end{lem}

For a graph $F$, let $\chi(F)$ denote its chromatic number. The following result is a direct corollary of Lemma \ref{spec-stability}.

\begin{cor} \label{var-spec-stability}
Let $\mathcal{F}$ be a finite family of graphs with $\min_{F\in \mathcal{F}}\chi(F)=r+1\geq3$. For every $\epsilon>0$, there exist $\delta>0$ and $n_{0}$ such that if $G$ is an $\mathcal{F}$-free graph of order $n\geq n_{0}$ with $\rho(G)>(\frac{r-1}{r}-\delta)n$, then $G$ can be
obtained from $T(n,r)$ by adding and deleting at most $\epsilon n^{2}$ edges.
\end{cor}

The following result is taken from \cite{Z1}.

\begin{lem}{\rm (\cite{Z1})}\label{mini degree}
Let $\mathcal{F}$ be a finite family of graphs with $\min_{F\in \mathcal{F}}\chi(F)=r+1\geq3$. For every $\theta>0$, there exists $n_{0}$ such that if $G\in {\rm SPEX}(n,\mathcal{F})$ with $n\geq n_{0}$, then $G$ is connected and $\delta(G)>(\frac{r-1}{r}-\theta)n$.
\end{lem}

The following lemma is taken from \cite{CFTZ}.

\begin{lem}{\rm (\cite{CFTZ})}\label{sets inequality}
Let $A_{1},A_{2},...,A_{\ell}$ be $\ell\geq2$ finite subsets of $A$. Then
$$|\cap_{1\leq i\leq \ell}A_{i}|\geq(\sum_{1\leq i\leq\ell}|A_{i}|)-(\ell-1)|\cup_{1\leq i\leq \ell}A_{i}|.$$
\end{lem}

The following lemma is taken from \cite{T}.

\begin{lem}{\rm (\cite{T})}\label{loop}
Let $H_{1}$ be a graph on $n_{0}$ vertices with maximum degree $d$ and $H_{2}$ be a
graph on $n-n_{0}$ vertices with maximum degree $d'$. $H_{1}$ and $H_{2}$ may have loops or multiple edges, where
loops add $1$ to the degree. Let $H=H_{1}\vee H_{2}$. Define
\begin{center}B =
$\left(\begin{array}{cc}
d&n-n_{0}\\
n_{0}&d'
\end{array}\right)$.
   \end{center}
Then $\rho(H)\leq\rho(B)=\frac{d+d'+\sqrt{(d-d')^{2}+4n_{0}(n-n_{0})}}{2}$.
\end{lem}

Furthermore, the equality in Lemma \ref{loop} can not hold, if either $H_{1}$ or $H_{2}$ is not regular.

\section{A general result for spectral extremal graphs}

A graph $F$ is called {\em vertex-critical}, if $\chi(F-u)=\chi(F)-1$ for some vertex $u\in V(F)$. Let $\mathcal{F}$ be a finite family of graphs with $\min_{F\in \mathcal{F}}\chi(F)=r+1\geq3$. $\mathcal{F}$ is called vertex-critical, if there is some $F_{0}\subseteq\mathcal{F}$ such that $F_{0}$ is vertex-critical with $\chi(F_{0})=r+1$.

\begin{theorem}\label{skeleton}
Assume that  $\mathcal{F}$ is a finite and vertex-critical graph family, where $\min_{F\in \mathcal{F}}\chi(F)=r+1\geq3$. Set $t=\max_{F\in\mathcal{F}}|F|$. For any small $\theta>0$, when $n$ is sufficiently large, each graph $G\in {\rm SPEX}(n,\mathcal{F})$  has the following conclusions.\\
$(i)$. There is a partition $V(G)=\cup_{1\leq i\leq r}V_{i}$ such that $||V_{i}|-\frac{n}{r}|<\theta n$ for any  $1\leq i\leq r$, and $d_{V_{i}}(v)< t$ and $d_{G}(v)\geq n-|V_{i}|$ for any $v\in V_{i}$.\\
$(ii)$. Let $\mathbf{x}=(x_{v})$ be a Perron vector of $G$ with the largest entry $1$. Then $x_{v}>1-\theta$ for any $v\in V(G)$.
\end{theorem}

\f{\bf Proof:} Let $G\in {\rm SPEX}(n,\mathcal{F})$ for large $n$, and let $\mathbf{x}=(x_{v})$ be a Perron vector of $G$ with the largest entry $1$. Since $\mathcal{F}$ is vertex-critical, there is some $F_{0}\subseteq\mathcal{F}$ such that $\chi(F_{0}-u)=r-1$ for some vertex $u\in V(F_{0})$. The following Claim 1 holds directly from Lemma \ref{mini degree}.

\medskip

\f{\bf Claim 1.}
 $G$ is connected and $\delta(G)>(\frac{r-1}{r}-\theta)n$.

\medskip

\f{\bf Claim 2.} $\rho(G)\geq\frac{r-1}{r}n-\frac{r}{4n}$.

\medskip

\f{\bf Proof of Claim 2.} Clearly, $T(n,r)$ is $\mathcal{F}$-free, since $\min_{F\in \mathcal{F}}\chi(F)=r+1$. Let $\overline{d}$ denote the average degree of $T(n,r)$. As is well known, $\overline{d}\geq\frac{r-1}{r}n-\frac{r}{4n}$.  It follows that $\rho(T(n,r))\geq\frac{r-1}{r}n-\frac{r}{4n}$. Then
$\rho(G)\geq\rho(T(n,r))\geq\frac{r-1}{r}n-\frac{r}{4n}$ as $G\in {\rm SPEX}(n,\mathcal{F})$. \hfill$\Box$

\medskip

\f{\bf Claim 3.} There is a partition $V(G)=\cup_{1\leq i\leq r}V_{i}$ such that $\sum_{1\leq i\leq r}e(G[V_{i}])$ is minimum. Moreover, $\sum_{1\leq i\leq r}e(G[V_{i}])<\theta^{3}n^{2}$ and $||V_{i}|-\frac{n}{r}|<\theta n$ for any $1\leq i\leq r$.

\medskip

\f{\bf Proof of Claim 3.} Since $G$ is $\mathcal{F}$-free and $\rho(G)\geq\frac{r-1}{r}n-\frac{r}{4n}$, by Corollary \ref{var-spec-stability} (letting $\epsilon=\frac{1}{2}\theta^{3}$), $G$ can be obtained from $T(n,r)$ by deleting and adding at most $\frac{1}{2}\theta^{3}n^{2}$ edges for large $n$. It follows that $e(G)>(\frac{r-1}{2r}-\theta^{3})n^{2}$. Moreover, there is a (balanced) partition $V(G)=\cup_{1\leq i\leq r}U_{i}$ such that $\sum_{1\leq i\leq r}e(G[U_{i}])<\theta^{3}n^{2}$. Now we select a partition $V(G)=\cup_{1\leq i\leq r}V_{i}$ such that $\sum_{1\leq i\leq r}e(G[V_{i}])$ is minimum. Then
 $$\sum_{1\leq i\leq r}e(G[V_{i}])\leq\sum_{1\leq i\leq r}e(G[U_{i}])<\theta^{3}n^{2}.$$

Let $a=\max_{1\leq i\leq r}||V_{i}|-\frac{n}{r}|$. Without loss of generality, assume that $a=||V_{1}|-\frac{n}{r}|$. Using the Cauchy-Schwarz inequality, we obtain that
$$2\sum_{2\leq i<j\leq r}|V_{i}||V_{j}|=(\sum_{2\leq i\leq r}|V_{i}|)^{2}-\sum_{2\leq i\leq r}|V_{i}|^{2}\leq\frac{r-2}{r-1}(\sum_{2\leq i\leq r}|V_{i}|)^{2}=\frac{r-2}{r-1}(n-|V_{1}|)^{2}.$$
Thus,
\begin{equation}
\begin{aligned}
e(G)&\leq(\sum_{1\leq i<j\leq r}|V_{i}||V_{j}|)+(\sum_{1\leq i\leq r}e(G[V_{i}]))\\
&\leq|V_{1}|(n-|V_{1}|)+(\sum_{2\leq i<j\leq r}|V_{i}||V_{j}|)+\theta^{3}n^{2}\\
&\leq|V_{1}|(n-|V_{1}|)+\frac{r-2}{2(r-1)}(n-|V_{1}|)^{2}+\theta^{3}n^{2}\\
&=\frac{r-1}{2r}n^{2}-\frac{r}{2(r-1)}a^{2}+\theta^{3}n^{2}.
\end{aligned}\notag
\end{equation}
Recall that $e(G)>(\frac{r-1}{2r}-\theta^{3})n^{2}$. It follows that $a\leq\sqrt{\frac{4(r-1)}{r}\theta^{3}n^{2}}<\theta n$ (requiring $\theta<\frac{r}{4(r-1)}$).
 This finishes the proof of Claim 3. \hfill$\Box$

\medskip

For $1\leq i\leq r$, let $W_{i}=\left\{v\in V_{i}~|~d_{V_{i}}(v)\geq\theta n\right\}$,  and let $W=\cup_{1\leq i\leq r}W_{i}$.

\medskip

\f{\bf Claim 4.}
$|W|< 2\theta^{2}n$.

\medskip

\f{\bf Proof of Claim 4.}  Since $\sum_{1\leq i\leq r}e(G[V_{i}])<\theta^{3}n^{2}$, and
$$\sum_{1\leq i\leq r}e(G[V_{i}])=\sum_{1\leq i\leq r}\frac{1}{2}\sum_{v\in V_{i}}d_{V_{i}}(v)\geq\sum_{1\leq i\leq r}\frac{1}{2}\sum_{v\in W_{i}}d_{V_{i}}(v)\geq\frac{1}{2}\sum_{1\leq i\leq r}\theta n|W_{i}|=\frac{1}{2}\theta n|W|,$$
we have $|W|<2\theta^{2}n$.
This finishes the proof of Claim 4.
\hfill$\Box$

\medskip

For any $1\leq i\leq r$, let $\overline{V}_{i}=V_{i}-W$.

\medskip

\f{\bf Claim 5.}
Let $1\leq \ell\leq r$ be fixed. For $i_{0}\neq \ell$, assume that $u_{0}\in W_{i_{0}}$ and  $u_{1},u_{2},...,u_{rt}\in \cup_{1\leq i\neq \ell\leq r}\overline{V}_{i}$. Then there are $t$ vertices in $\overline{V}_{\ell}$ which are adjacent to all the vertices $u_{0},u_{1},u_{2},...,u_{rt}$ in $G$.

\medskip

\f{\bf Proof of Claim 5.} Recall that $\frac{n}{r}-\theta n\leq|V_{s}|\leq\frac{n}{r}+\theta n$ for any $1\leq s\leq r$ by Claim 3. By Claim 4, we have $\frac{n}{r}-2\theta n\leq|\overline{V}_{s}|\leq\frac{n}{r}+\theta n$. By Claim 1, we have $\delta(G)>(\frac{r-1}{r}-\theta)n$. Since $\sum_{1\leq i\leq r}e(G[V_{i}])$ is minimum, we have $d_{V_{\ell}}(u_{0})\geq d_{V_{i_{0}}}(u_{0})$  as $u_{0}\in V_{i_{0}}$.  Otherwise, we will obtain a contradiction by moving $u_{0}$ from $V_{i_{0}}$ to $V_{\ell}$. It follows that
$$d_{V_{\ell}}(u_{0})\geq\frac{d_{V_{\ell}}(u_{0})+d_{V_{i_{0}}}(u_{0})}{2}\geq\frac{d_{G}(u_{0})-\sum_{1\leq j\neq\ell,i_{0}\leq r}|V_{j}|}{2}\geq(\frac{1}{2r}-(r-1)\theta)n.$$
Then, using Claim 4,
 $$d_{\overline{V}_{\ell}}(u_{0})\geq d_{V_{\ell}}(u_{0})-|W|\geq(\frac{1}{2r}-r\theta)n.$$

 For any $1\leq i\leq rt$, assume that $u_{i}\in \overline{V}_{j_{i}}$, where $j_{i}\neq\ell$. Then $d_{V_{j_{i}}}(u_{i})\leq\theta n$ as $u_{i}\notin W$.
Hence $$d_{V_{\ell}}(u_{i})\geq d_{G}(u_{i})-d_{V_{j_{i}}}(u_{i})-\sum_{1\leq s\neq\ell,j_{i}\leq r}|V_{s}|\geq(\frac{1}{r}-r\theta)n.$$
Thus,
 $$d_{\overline{V}_{\ell}}(u_{i})\geq d_{V_{\ell}}(u_{i})-|W|\geq(\frac{1}{r}-2r\theta)n.$$
By Lemma \ref{sets inequality}, we have
$$|N_{\overline{V}_{\ell}}(u_{0})\cap(\cap_{1\leq i\leq rt}N_{\overline{V}_{\ell}}(u_{i}))|\geq|N_{\overline{V}_{\ell}}(u_{0})|+(\sum_{1\leq i\leq rt}|N_{\overline{V}_{\ell}}(u_{i})|)-rt|V_{\ell}|\geq(\frac{1}{2r}-rt(r+2)\theta)n\geq t.$$
Thus, there are $t$ vertices in $\overline{V}_{\ell}$ which are adjacent to all the vertices $u_{0},u_{1},u_{2},...,u_{rt}$.
 This finishes the proof of Claim 5. \hfill$\Box$

\medskip

\medskip

\f{\bf Claim 6.}
$W=\emptyset$. Moreover, $d_{V_{i}}(v)<t$ for any $v\in V_{i}$ and $1\leq i\leq r$.

\medskip

\f{\bf Proof of Claim 6.} Suppose that $W\neq\emptyset$. Let $w\in W_{1}$ without loss of generality. Since $d_{V_{1}}(w)\geq \theta n$ and $|W|\leq2\theta^{2}n$, $w$ has $\theta n-2\theta^{2}n\geq t$ neighbors in $\overline{V}_{1}$, say $u_{1},u_{2},...,u_{t}$. By Claim 5, $w,u_{1},u_{2},...,u_{t}$ have $t$ common neighbors in $\overline{V}_{2}$, say $u_{t+1},u_{t+2},...,u_{2t}$. Repeat the process using Claim 5. We can obtain a copy of $F_{0}$ in $G$, a contradiction. Thus $W=\emptyset$. Then $\overline{V}_{i}=V_{i}$ for any $1\leq i\leq r$.

If $d_{V_{i_{1}}}(v)\geq t$ for some $v\in V_{i_{1}}$ and $1\leq i_{1}\leq r$, then $v$ has $t$ neighbors in $V_{i_{1}}$, say $u_{1},u_{2},...,u_{t}$. By Claim 5, $v,u_{1},u_{2},...,u_{t}$ have $t$ common neighbors in $V_{i_{2}}$ with $i_{2}\neq i_{1}$, say $u_{t+1},u_{t+2},...,u_{2t}$. Again by Claim 5, $v,u_{1},u_{2},...,u_{2t}$ have $t$ common neighbors in $V_{i_{3}}$ with $i_{3}\neq i_{1},i_{2}$, say $u_{2t+1},u_{2t+2},...,u_{3t}$. Repeat the process using Claim 5. We can obtain a copy of $F_{0}$ in $G$, a contradiction. Hence $d_{V_{i}}(v)< t$ for any $v\in V_{i}$ and $1\leq i\leq r$.
 This finishes the proof of Claim 6. \hfill$\Box$

\medskip

\f{\bf Claim 7.}
$x_{v}>1-\theta$ for any $v\in V(G)$.

\medskip

\f{\bf Proof of Claim 7.} Recall that $\mathbf{x}=(x_{v})$ has the largest entry $1$. Without loss of generality, assume that $x_{v^{*}}=1$ and $v^{*}\in V_{1}$. Note that $d_{V_{1}}(v^{*})<t$ by Claim 6.
Then
\begin{equation}
\begin{aligned}
\rho(G)x_{v^{*}}&=(\sum_{v\in N_{V_{1}}(v^{*})}x_{v})+(\sum_{v\in N_{\cup_{i\neq1}V_{i}}(v^{*})}x_{v})\\
&\leq t-1+\sum_{v\in N_{\cup_{i\neq1}V_{i}}(v^{*})}x_{v}\\
&\leq t-1+\sum_{v\in \cup_{i\neq1}V_{i}}x_{v}.
\end{aligned}\notag
\end{equation}
It follows that
$$\sum_{v\in \cup_{i\neq1}V_{i}}x_{v}\geq\rho(G)-t+1.$$

Suppose that $x_{u'}<\frac{\rho(G)-t+1}{1+\rho(G)}$ for some $u'\in V(G)$.  Then
$$x_{u'}+\sum_{v\in N_{G}(u')}x_{v}=(1+\rho(G))x_{u'}<\rho(G)-t+1\leq\sum_{v\in \cup_{i\neq1}V_{i}}x_{v}.$$
Let $G'$ be the graph obtained
from $G$ by deleting all the edges incident with $u'$ and adding all the edges between $u'$ and $(\cup_{i\neq1}V_{i})-\left\{u'\right\}$.
Then
$$\mathbf{x}^{T}(\rho(G')-\rho(G))\mathbf{x}\geq\mathbf{x}^{T}(A(G')-A(G))\mathbf{x}=2x_{u'}((\sum_{v\in \cup_{i\neq1}V_{i}}x_{v})-(x_{u'}+\sum_{v\in N_{G}(u')}x_{v}))>0,$$
implying that $\rho(G')>\rho(G)$. Now we show that $G'$ is $\mathcal{F}$-free. In fact,
if $F\subseteq G'$ for some $F\in \mathcal{F}$, then $u'\in V(F)$. Note that the neighbors of $u'$ in $F$ are all contained in $\cup_{i\neq1}V_{i}$, say $u_{1},u_{2},...,u_{t_{0}}$ with $1\leq t_{0}\leq t$. By Claim 5, $u_{1},u_{2},...,u_{t_{0}}$ have at least one common neighbor $u''$ in $V_{1}$ such that $u''\notin V(F)$. Let $F'$ be obtained from $F$ by deleting $u'$ and adding $u''$. Clearly, $F\subseteq F'$ and $F'\subseteq G$. That is, $F\subseteq G$, a contradiction. Hence, $G'$ is $\mathcal{F}$-free. But since $\rho(G')>\rho(G)$, it contradicts that $G\in {\rm SPEX}(n,\mathcal{F})$. Thus, we must have $x_{v}\geq\frac{\rho(G)-t+1}{1+\rho(G)}$ for any $v\in V(G)$. By Claim 2, $\rho(G)>\frac{n}{3}$ for large $n$. It follows that $x_{v}>1-\theta$ for any $v\in V(G)$.
This completes the proof of Claim 7. \hfill$\Box$

\medskip

To complete the proof, it remains to show $d_{G}(v)\geq n-|V_{i}|$ for any $v\in V_{i}$ and $1\leq i\leq r$. Without loss of generality, suppose that $d_{G}(v)< n-|V_{1}|$ for some $v\in V_{1}$. Set $d=d_{V_{1}}(v)$. Then $v$ has at least $d+1$ non-neighbors in $\cup_{i\neq1}V_{i}$, say $w_{1},w_{2},...,w_{d+1}$. Let $G''$ be the graph obtained from $G$ by deleting  the $d$ edges incident with $v$ and inside $V_{1}$, and adding $d_{V_{1}}(v)+1$ non-edges between $v$ and $\cup_{i\neq1}V_{i}$. Since $1-\theta\leq x_{w}\leq1$ for any $w\in V(G)$ by Claim 7, it is easy to see
\begin{equation}
\begin{aligned}
&\mathbf{x}^{T}(\rho(G'')-\rho(G))\mathbf{x}\\
&\geq\mathbf{x}^{T}(A(G'')-A(G))\mathbf{x}\\
&=2x_{v}((\sum_{1\leq j\leq d+1}x_{w_{j}})-(\sum_{w\in N_{V_{1}}(v)}x_{w}))\\
&\geq2(1-\theta)((d+1)(1-\theta)-d)\\
&>0.
\end{aligned}\notag
\end{equation}
It follows that $\rho(G'')>\rho(G)$. Similar to Claim 7, we can show that $G''$ is $\mathcal{F}$-free.
 But it contradicts that $G\in {\rm SPEX}(n,\mathcal{F})$. Hence, we must have $d_{G}(v)\geq n-|V_{i}|$ for any $v\in V_{i}$ and $1\leq i\leq r$.
This completes the proof. \hfill$\Box$

\section{Spectral radius and walks}

To prove Theorem \ref{main1}, we need a result in \cite{Z}. 
Let $G$ be a graph. For an integer $\ell\geq1$, $v_{0}v_{1}\cdots v_{\ell}$ is called a {\em walk} of length $\ell$ in $G$, if $v_{i}\sim v_{i+1}$ for any $0\leq i\leq \ell-1$. The vertex $v_{0}$ is called the starting vertex. For any $u\in V(G)$, let $w^{\ell}_{G}(u)$ be the number of walks of length $\ell$ starting at $u$. Let $W^{\ell}(G)=\sum_{v\in V(G)}w^{\ell}_{G}(u)$.
 For any integers $\ell\geq2$ and $1\leq i\leq \ell-1$, the following formula (by considering the $(i+1)$-th vertex in a walk of length $\ell$) will be used:
 $$W^{\ell}(G)=\sum_{u\in V(G)}w^{i}_{G}(u)w^{\ell-i}_{G}(u).$$

For two graphs $G_{1}$ and $G_{2}$, we say $G_{1}\succ G_{2}$, if there is an integer $\ell\geq1$ such that $W^{\ell}(G_{1})>W^{\ell}(G_{2})$ and $W^{i}(G_{1})=W^{i}(G_{2})$ for any $1\leq i\leq \ell-1$;  $G_{1}\equiv G_{2}$, if $W^{i}(G_{1})=W^{i}(G_{2})$ for any $i\geq1$; $G_{1}\prec G_{2}$ , if $G_{2}\succ G_{1}$.

For a family of graphs $\mathcal{G}$,
let $${\rm EX}^{1}(\mathcal{G})=\left\{G\in\mathcal{G}~|~W^{1}(G)\geq W^{1}(G')~ {\rm for~any}~G'\in\mathcal{G}\right\},$$
 and
$${\rm EX}^{\ell}(\mathcal{G})=\left\{G\in{\rm EX}^{\ell-1}(\mathcal{G})~|~W^{\ell}(G)\geq W^{\ell}(G')~ {\rm for~any}~G'\in{\rm EX}^{\ell-1}(\mathcal{G})\right\}$$
 for any $\ell\geq2$. By definition, ${\rm EX}^{i+1}(\mathcal{G})\subseteq {\rm EX}^{i}(\mathcal{G})$ for any $i\geq1$. Let ${\rm EX}^{\infty}(\mathcal{G})=\cap_{1\leq i\leq\infty}{\rm EX}^{i}(\mathcal{G})$.

\medskip

The following result is taken from \cite{Z}.

\begin{theorem}{\rm (\cite{Z})}\label{one-set}
Let $G$ be a connected graph of order $n$, and let $S$ be a subset of $V(G)$ with $1\leq|S|<n$. Assume that $T$ is a set of some isolated vertices of $G-S$, such that each vertex in $T$ is adjacent to each vertex in $S$ in $G$. Let $H_{1}$ and $H_{2}$ be two graphs with vertex set $T$. For any $1\leq i\leq2$, let $G_{i}$ be the graph obtained from $G$ by embedding the edges of $H_{i}$ into $T$. When $\rho(G)$ is sufficiently large (compared with $|T|$), we have the following conclusions.\\
$(i)$ If $H_{1}\equiv H_{2}$, then $\rho(G_{1})=\rho(G_{2})$.\\
$(ii)$ If $H_{1}\succ H_{2}$, then $\rho(G_{1})>\rho(G_{2})$.\\
$(iii)$ If $H_{1}\prec H_{2}$, then $\rho(G_{1})<\rho(G_{2})$.
\end{theorem}

The following lemma is taken from \cite{IS}.

\begin{lem}{\rm (\cite{IS})}\label{2-degree}
Let $\mathcal{M}_{n,m}$ be the set of all the graphs of order $n$ with $m$ edges, where $m\geq1$ and $n\geq m+2$. For all graphs $G\in\mathcal{M}_{n,m}$, $\sum_{v\in V(G)}d^{2}_{G}(u)$ is maximized  when $G\in\left\{K_{1,3}\cup \overline{K_{n-4}},K_{3}\cup \overline{K_{n-3}}\right\}$ for $m=3$, and $G=\left\{K_{1,m}\cup \overline{K_{n-1-m}}\right\}$ otherwise.
\end{lem}

\begin{lem}\label{bounded order}
For odd integers $k\geq3$ and $n>2k+1$, let $G$ be a connected nearly $k$-regular graph of order $n$. Let $u$ be the unique vertex with degree $k-1$. If the neighbors of $u$ induce a clique in $G$, then $G$ contains a $P_{2k+3}$.
\end{lem}

\f{\bf Proof:} Since $n>2k+1$ is odd, we have $n\geq2k+3$. Let $u_{1},u_{2},...,u_{k-1}$ be the neighbors of $u$. By assumption, $G[\left\{u_{1},u_{2},...,u_{k-1}\right\}]$ is a complete graph. Let $H=G-\left\{u,u_{1},u_{2},...,u_{k-1}\right\}$. Recall that all the vertices except $u$ have degree $k$ in $G$. Hence, $u_{i}$ has exactly one neighbor in $V(H)$ for each  $1\leq i\leq k-1$. We will prove the lemma by several cases.

\medskip

\f{\bf Case 1.} There is exactly one vertex in $V(H)$, say $v$, which has neighbors in $\left\{u_{1},u_{2},...,u_{k-1}\right\}$.

\medskip

In this case, $v$ is adjacent to all the vertices in $\left\{u_{1},u_{2},...,u_{k-1}\right\}$. Thus, $v$ has exactly one neighbor in $H$, say $w$. Let $H'=H-v$. Clearly, $H'$ is a connected nearly $k$-regular graph with the unique vertex $w$ of degree $k-1$. Since $n\geq2k+3$, we have $|H'|\geq k+2$. Let $Q=ww_{1}w_{2}\cdots w_{\ell}$ be a longest path starting at $w$ in $H'$. We shall prove $\ell\geq k+1$. Suppose that $\ell\leq k$. Since $Q$ is longest, all the neighbors of $w_{\ell}$ are in $\left\{w,w_{1},w_{2},...,w_{\ell-1}\right\}$. This implies that $\ell\geq k$, and thus $\ell=k$. Moreover, $w_{k}(=w_{\ell})$ is adjacent to all other vertices in $\left\{w,w_{1},w_{2},\cdots, w_{k}\right\}$. But then, $ww_{k}w_{1}w_{2}\cdots w_{k-1}$ is also a longest path starting at $w$ in $H'$. Similarly, $w_{k-1}$ is adjacent to all other vertices in $\left\{w,w_{1},w_{2},\cdots, w_{k}\right\}$.  But then, $ww_{k-1}w_{k}w_{1}w_{2}\cdots w_{k-2}$ is also a longest path starting at $w$ in $H'$. Repeating this process, we can obtain that $w_{1},w_{2},\cdots, w_{k}$ are all adjacent to $w$.  This contradicts the fact that $w$ has degree $k-1$ in $H'$. Hence, $\ell\geq k+1$ is proved. Then $uu_{1}u_{2}\cdots u_{k-1}vww_{1}w_{2}\cdots w_{\ell}$ is a path of order at least $2k+3$ in $G$, as desired.

 \medskip
 
\f{\bf Case 2.} There are at least two vertices in $V(H)$, which have neighbors in $\left\{u_{1},u_{2},...,u_{k-1}\right\}$. 

\medskip

  Without loss of generality, assume that $x\in V(H)$ is adjacent to $u_{1}$, and $y\in V(H)$ is adjacent to $u_{k-1}$, where $x\neq y$. (Note that each $u_{i}$ has exactly one neighbor in $V(H)$.) Recall that $|H|\geq k+3$ as $n\geq2k+3$.

\medskip

\f{\bf Subcase 2.1.} $H$ is connected.

\medskip

If $H$ contains a Hamilton cycle, then clearly, $G$ contains a Hamilton path, as desired. Thus, we can assume that $H$ contains no Hamilton cycles.
 Let $Q=v_{1}v_{2}\cdots v_{\ell}$ be a longest path in $H$. Then each neighbor in $H$ of $v_{1}$ and $v_{\ell}$ is contained in $Q$. Let $d_{1}$ and $d_{\ell}$ be the numbers of neighbors in $\left\{u_{1},u_{2},...,u_{k-1}\right\}$ of $v_{1}$ and $v_{\ell}$, respectively. We can require that $d_{1}+d_{\ell}\leq k-2$. In fact, it is clear that $d_{1}+d_{\ell}\leq k-1$. Suppose that $d_{1}+d_{\ell}= k-1$. Without loss of generality, we can assume that $d_{\ell}\leq d_{1}$, implying that $d_{\ell}\leq\frac{k-1}{2}$. Hence, $v_{\ell}$ has at least $\frac{k+1}{2}\geq2$ neighbors in $H$. So, there is a $v_{i}$ adjacent to $v_{\ell}$, where $1\leq i< \ell-1$. Then $v_{i+1}v_{i+2}\cdots v_{\ell}v_{i}v_{i-1}\cdots v_{1}$ is a longest path in $H$. We can use $v_{i+1}$ instead of $v_{\ell}$ if necessary. Hence, we can assume that $d_{1}+d_{\ell}\leq k-2$.
 
 Set $U=\left\{i-1~|~v_{1}\sim v_{i},2\leq i\leq \ell\right\}$ and $W=\left\{j~|~v_{\ell}\sim v_{j},1\leq j\leq \ell-1\right\}$. Clearly, $|U|=k-d_{1}$ and $|W|=k-d_{\ell}$. Now we prove that $ U\cap W=\emptyset$. Otherwise,
suppose that $i_{0}\in U\cap W$. Then $v_{1}\sim v_{i_{0}+1}$ and  $v_{\ell}\sim v_{i_{0}}$.  It follows that $v_{1},v_{2},...,v_{i_{0}}v_{\ell},v_{\ell-1},...,v_{i_{0}+1}v_{1}$ is a cycle of order $\ell$ in $H$. Since $H$ is a connected graph without Hamilton cycles, there is another vertex in $H$ adjacent to some vertex in this cycle. Then we can obtain a path of order $\ell+1$ in $H$. This is impossible as $Q$ is a longest path in $H$. Hence, $ U\cap W=\emptyset$ is proved. Note that $U,W\subseteq\left\{1,2,...,\ell-1\right\}$. It follows that $\ell\geq1+|U|+|W|=2k+1-(d_{1}+d_{\ell})$. 

\medskip

\f{\bf Subcase 2.1.1.} $d_{1}+d_{\ell}\geq1$.

\medskip

Without loss of generality, assume that $d_{1}\geq1$. Then $v_{1}$ has a neighbor in $\left\{u_{1},u_{2},...,u_{k-1}\right\}$, say $u_{k-1}$. Recall that $d_{1}+d_{\ell}\leq k-2$. Thus, $\ell\geq2k+1-(d_{1}+d_{\ell})\geq k+3$. Hence, $uu_{1}u_{2}\cdots u_{k-1}v_{1}v_{2}\cdots v_{\ell}$ is a path of order at least $2k+3$, as desired.

\medskip

\f{\bf Subcase 2.1.2.} $d_{1}=d_{\ell}=0$.

\medskip

In this case, $\ell\geq2k+1$. Recall that $x\in V(H)$ is adjacent to $u_{1}$, and $y\in V(H)$ is adjacent to $u_{k-1}$. 

\medskip

\f{\bf Subcase 2.1.2.1.} There is at least one vertex of $x$ and $y$, which is in $V(H)-V(Q)$. 

\medskip

Without loss of generality, assume that $x$ is in $V(H)-V(Q)$. First consider $y\neq v_{k+1}$. Let $Q_{y}$ be a shortest path in $H$ from $y$ to $Q$. ($Q_{y}=y$ if $y\in V(Q)$.) We can require that $Q_{y}$ is not through $x$, since we can change $x$ and $y$ if necessary. Let $v_{j}$ be the other end vertex of $Q_{y}$, where $1\leq j\leq \ell$. Denote $yQ_{y}v_{j}=Q_{y}$. If $j\leq k+1$, then (noting that $y\neq v_{k+1}$), $xu_{1}uu_{2}u_{3}\cdots u_{k-1}yQ_{y}v_{j}v_{j+1}\cdots v_{\ell}$ is a path of order at least $2k+3$, as desired. If $j>k+1$, then $xu_{1}uu_{2}u_{3}\cdots u_{k-1}yQ_{y}v_{j}v_{j-1}\cdots v_{1}$ is a path of order at least $2k+3$, as desired. 

It remains to consider $y=v_{k+1}$. Since $x$ has at most $k-2$ neighbors in $\left\{u_{1},u_{2},...,u_{k-1}\right\}$, $x$  has at least 2 neighbors in $V(H)$. Let $x'$ be a neighbor in $V(H)$ of $x$, such that $x'\neq v_{k+1}$. Clearly, at least one of $v_{1}v_{2}\cdots v_{k+1}$ and $v_{k+1}v_{k+2}\cdots v_{\ell}$ is not through $x'$. Without loss of generality, assume that $v_{k+1}v_{k+2}\cdots v_{\ell}$ is not through $x'$. Then, 
$$x'xu_{1}uu_{2}u_{3}\cdots u_{k-1}v_{k+1}v_{k+2}\cdots v_{\ell}$$
 is a path of order at least $2k+3$, as desired.

\medskip

\f{\bf Subcase 2.1.2.2.} Both $x$ and $y$ are in $V(Q)$.

\medskip

First consider the case that $x$ or $y$ is $v_{i}$ for some $i\leq k-1$ or $i\geq k+3$. Without loss of generality (it is very similar for other cases), we can assume that $x=v_{i}$ with $i\leq k-1$. Then, $uu_{k-1}u_{k-2}\cdots u_{1}v_{i}v_{i+1}v_{i+2}\cdots v_{\ell}$ is a path of order at least $2k+3$, as desired.

It remains that $x,y\in\left\{v_{k},v_{k+1},v_{k+2}\right\}$. Without loss of generality (it is very similar for other cases), we can assume that $x=v_{k}$ and $y=v_{k+2}$. Then 
$$v_{1}v_{2}\cdots v_{k}u_{1}uu_{2}u_{3}\cdots u_{k-1}v_{k+2}v_{k+3}\cdots v_{\ell}$$
 is a path of order at least $3k\geq2k+3$, as desired.

\medskip

\f{\bf Subcase 2.2.} $H$ is not connected.

\medskip

Let $H^{1}$ and $H^{2}$ be two components of $H$. Since $G$ is connected, each component of $H$ has at least one neighbor in $\left\{u_{1},u_{2},...,u_{k-1}\right\}$. We can assume that $x\in V(H^{1})$ and $y\in V(H^{2})$ without loss of generality. Let $Q_{1}=xx_{1}x_{2}\cdots x_{a}$ be a longest path starting at $x$ in $H^{1}$, and let $Q_{2}=yy_{1}y_{2}\cdots y_{b}$ be a longest path starting at $y$ in $H^{2}$. Let $d_{a}$ and $d_{b}$ be the numbers of  neighbors in $\left\{u_{1},u_{2},...,u_{k-1}\right\}$ of $x_{a}$ and $y_{b}$, respectively. Then $d_{a}+d_{b}\leq k-1-2=k-3$. Since $Q_{1}$ and $Q_{2}$ are longest, each neighbor of $x_{a}$ or $y_{b}$ is in $Q_{1}$ or $Q_{2}$. It follows that $a\geq k-d_{a}$ and $b\geq k-d_{b}$. Thus, $a+b\geq 2k-(d_{a}+d_{b})\geq k+3$. Then, 
$x_{a}x_{a-1}\cdots xu_{1}uu_{2}u_{3}\cdots u_{k-1}yy_{1}y_{2}\cdots y_{b}$ is a path of order at least $a+1+k+b+1\geq2k+5$, as desired.
This completes the proof. \hfill$\Box$

\medskip

 For odd integer $k\geq3$, let $Q^{*}_{k}$ be the graph obtained from $K_{1}\vee K_{k-1}$ and $K_{\frac{k-1}{2}}\vee K_{\frac{k+1}{2}}$ by adding a single vertex $w$, and connecting $w$ to all vertices in the part $V(K_{\frac{k+1}{2}})$ and to $\frac{k-1}{2}$ vertices in $V(K_{k-1})$, and then connecting the remained $\frac{k-1}{2}$ vertices in $V(K_{k-1})$ to the vertices in the part $V(K_{\frac{k-1}{2}})$ by a matching of $\frac{k-1}{2}$ edges. Clearly, $Q^{*}_{k}$ is a nearly $k$-regular graph of order $2k+1$. 

For $k=7$, let $Q^{**}$ be the graph obtained from $K_{1}\vee K_{6}$ and $\overline{K}_{3}\vee K_{5}$ by connecting each vertex in the part $V(\overline{K}_{3})$ to exactly two vertices in the part $V(K_{6})$. Clearly, $Q^{**}$ is a nearly $7$-regular graph of order $15$.

For odd integers $k\geq3$ and $n\geq4k+3$, let $\mathcal{G}_{n,k}$ be the set of nearly $k$-regular $P_{2k+3}$-free graphs of order $n$.
Let $\mathcal{V}_{n,k}$ be the set of nearly $k$-regular $P_{2k+3}$-free graphs $G$ of order $n$, where $G$ has a component $Q^{*}_{k}$ for $k\neq7$ and has a component $Q^{**}$ for $k=7$. Clearly, $\mathcal{V}_{n,k}\subseteq \mathcal{G}_{n,k}$, and all the graphs in $\mathcal{V}_{n,k}$ have the same number of walks of length $\ell$ for any $\ell\geq1$.

The following fact will be used in the rest of this paper.

\medskip

\f{\bf Fact 1.} For any two integers $k\geq1$ and $m\geq k+1$, there is a $k$-regular graph of order $m$ if and only if $km$ is even.

\medskip

\begin{lem}\label{walk lemma}
Let $k\geq3$ and $n\geq4k+3$ be odd integers. Then ${\rm EX}^{\infty}(\mathcal{G}_{n,k})=\mathcal{V}_{n,k}$.
\end{lem}

\f{\bf Proof:} Recall that there exists a $k$-regular graph of order $m$ if and only if $k m$ is even and $m\geq k+1$. Since $|Q^{*}_{k}|=2k+1$ for any $k\geq3$ and $|Q^{**}|=|Q^{*}_{7}|=15$, we have that $n-(2k+1)\geq2(k+1)$ can be partitioned into several even integers between $k+1$ and $2k$. Thus, $\mathcal{V}_{n,k}$
is not empty by Fact 1. Let $H\in\mathcal{V}_{n,k}$. Then the component of $H$ including the vertex with degree $k-1$ is $Q^{*}_{k}$ for $k\neq7$ and $Q^{**}$ for $k=7$.

Let $G\in {\rm EX}^{\infty}(\mathcal{G}_{n,k})$. Let $u$ be the vertex with degree $k-1$ in $G$, and let $Q$ be the component of $G$ including $u$. We will prove that $Q=Q^{*}_{k}$ for $k\neq7$ and $Q=Q^{**}$ for $k=7$. Then $G\in\mathcal{V}_{n,k}$, and the lemma follows. 

Set $q=|Q|$. Then $q$ is odd, since $Q$ is a nearly $k$-regular graph. For any $1\leq i\leq 2$, let $N_{i}$ be the set of vertices of $Q$ at distance $i$ from the vertex $u$. Let $N_{3}$ be the set of vertices of $Q$ at distance $\geq3$ from $u$. Then $|N_{1}|=k-1$. 
For any $v\in V(Q)$ and $1\leq i\leq3$, let $d_{i}(v)$ be the number of neighbors of $v$ in $N_{i}$. Set $e(N_{1},N_{2})=e_{G}(N_{1},N_{2})$.
Clearly, $$e(N_{1},N_{2})\geq\sum_{v\in N_{1}}d_{2}(v)|N_{1}|\geq|N_{1}|=k-1.$$ Moreover, $e(N_{1},N_{2})=k-1$ if and only if $Q[N_{1}]$ is a clique.

   Let $w^{i}(v)=w^{i}_{G}(v)$ for any $i\geq1$ and $v\in V(G)$. It is easy to check that:
 $$w^{2}(u)=k^{2}-k,w^{3}(u)=(k^{2}-1)(k-1);$$
for $v\in N_{1}$,
$$w^{2}(v)=k^{2}-1,w^{3}(v)=k^{2}-k+d_{1}(v)(k^{2}-1)+d_{2}(v)k^{2}=k^{3}-2k+1+d_{2}(v);$$
 for $v\in N_{2}$,
$$w^{2}(v)=k^{2},w^{3}(v)=d_{1}(v)(k^{2}-1)+(k-d_{1}(v))k^{2}=k^{3}-d_{1}(v);$$
for any $v\in N_{3}$ or $v\in V(G)-V(Q)$ and $1\leq i\leq3$,
$$w^{i}(v)=k^{i}.$$
Then by a calculation, we have that
$$W^{1}(G)=nk-1,$$
$$W^{2}(G)=(n-1)k^{2}+(k-1)^{2},$$
$$W^{3}(G)=\sum_{v\in V(G)}w^{1}(v)w^{2}(v)=nk^{3}-3k^{2}+2k,$$
$$W^{4}(G)=\sum_{v\in V(G)}w^{2}(v)w^{2}(v)=nk^{4}-4k^{3}+3k^{2}+k-1,$$
$$W^{5}(G)=\sum_{v\in V(G)}w^{2}(v)w^{3}(v)=(n-q)k^{5}+\sum_{v\in V(Q)}w^{2}(v)w^{3}(v),$$
$$=(n-k)k^{5}+k(k+1)(k-1)^{3}+(k^{2}-1)(k-1)(k^{3}-2k+1)-e(N_{1},N_{2})$$
and
$$W^{6}(G)=\sum_{v\in V(G)}w^{3}(v)w^{3}(v)=(n-q)k^{6}+\sum_{v\in V(Q)}w^{3}(v)w^{3}(v)$$
$$=(n-k)k^{6}+(k-1)^{2}(k^{2}-1)^{2}+(k-1)(k^{3}-2k+1)^{2}-(4k-2)e(N_{1},N_{2})
$$
$$+\sum_{v\in N_{1}}d_{2}^{2}(v)+\sum_{v\in N_{2}}d_{1}^{2}(v).$$

Using a similar discussion for $H$, we can obtain that $W^{i}(G)=W^{i}(H)$ for $1\leq i\leq4$. Since $G\in {\rm EX}^{\infty}(\mathcal{G}_{n,k})$, we have $W^{5}(G)\geq W^{5}(H)$. This requires that $e(N_{1},N_{2})\geq k-1$ is minimized. Then $e(N_{1},N_{2})=k-1$ and $W^{5}(G)= W^{5}(H)$. Moreover, $G[N_{1}]$ is a clique of order $k-1$. Also by $G\in {\rm EX}^{\infty}(\mathcal{G}_{n,k})$, we have $W^{6}(G)\geq W^{6}(H)$ as $W^{i}(G)= W^{i}(H)$ for $1\leq i\leq5$.
This requires that $\sum_{v\in N_{1}}d_{2}^{2}(v)+\sum_{v\in N_{2}}d_{1}^{2}(v)$ is maximized. Note that $d_{2}(v)=1$ for any $v\in N_{1}$. This requires that $\sum_{v\in N_{2}}d_{1}^{2}(v)$ is maximized, subject to $\sum_{v\in N_{2}}d_{1}(v)=k-1$.

By the above discussion, $N_{1}$ induces a clique in $Q$. Since $Q$ contains no $P_{2k+3}$, we have that $q\leq2k+1$ by Lemma \ref{bounded order}. Since $e(N_{1},N_{2})=k-1$, we must have $q=2k+1$. Let $Q'=\overline{Q-(\left\{u\right\}\cup N_{1})}$ (i.e., the complement of $Q-(\left\{u\right\}\cup N_{1})$). Clearly, $d_{Q'}(v)=d_{1}(v)$ for any $v\in N_{2}$. This requires that $\sum_{v\in V(Q')}d_{Q'}^{2}(v)$ is maximized, subject to $\sum_{v\in V(Q')}d_{Q'}(v)=k-1$. By Lemma \ref{2-degree}, we must have $Q'=K_{1,\frac{k-1}{2}}\cup\overline{K_{\frac{k+1}{2}}}$ for $k\geq3$, or $Q'=K_{3}\cup\overline{K_{5}}$ for $k=7$. Thus, $Q=Q^{*}_{k}$ for $k\geq3$ or $Q=Q^{**}$ for $k=7$. Both $Q^{*}_{7}$ and $Q^{**}$ are nearly $7$-regular graphs of order $15$. By a calculation, we can show $W^{7}(Q^{**})-W^{7}(Q^{*}_{7})=84>0$ (using formula $W^{7}(G)=\sum_{v\in V(G)}w^{3}(v)w^{4}(v)$). Hence, $Q=Q^{**}$ for $k=7$. Consequently, $Q=Q^{*}_{k}$ for $k\neq7$. This completes the proof. \hfill$\Box$

\medskip

For $k\geq1$, define $f(k)=\sum_{0\leq i\leq 2k+1}k^{i}$. Clearly, any connected $P_{2k+3}$-free graph $G$ of order $n$ with $\Delta(G)\leq k$, has diameter at most $2k+1$. Thus, $n\leq f(k)$.

\begin{theorem}\label{spec compare}
Assume that $k\geq3$ is odd and  $n\equiv 2(\mod4)$ is sufficiently large. Let $G$ be a graph obtained from the Tur\'{a}n graph $T(n,2)$ with parts $L$ and $R$, by embedding a graph from $\mathcal{G}_{\frac{n}{2},k}$ into $G[L]$. Then $\rho(G)<\frac{k+\sqrt{k^{2}+n^{2}-4}}{2}$.
\end{theorem}

\f{\bf Proof:} Since $G[L]\in \mathcal{G}_{\frac{n}{2},k}$, each component of $G[L]$ has order at most $f(k)$.
 Let $Q_{1}$ be the component of $G[L]$ including the unique vertex with degree $k-1$.
Since each component of $G[L]$ has order between $k+1$ and $f(k)$, we can choose several other components $Q_{2},Q_{3},...,Q_{t}$, such that $4k+3\leq|\cup_{1\leq i\leq t}Q_{i}|\leq4k+3+f(k)$. Let $Q=\cup_{1\leq i\leq t}Q_{i}$. Since $Q\in \mathcal{G}_{|Q|,k}$ and $n$ is sufficiently large with respect to $|Q|$, by Theorem \ref{one-set}, $\rho(G)$ is maximized when $Q\in{\rm EX}^{\infty}(\mathcal{G}_{|Q|,k})$. That is, $Q\in\mathcal{V}_{|Q|,k}$ by Lemma \ref{walk lemma}. And then $G[L]\in \mathcal{V}_{\frac{n}{2},k}$.  To prove the theorem, we can assume that $G[L]\in \mathcal{V}_{\frac{n}{2},k}$.

Let $g(x)=x^{2}-kx-(\frac{n^{2}}{4}-1)$. Clearly, $x=\frac{k+\sqrt{k^{2}+n^{2}-4}}{2}$ is the largest root of $g(x)=0$. 

For $k=7$, $G[L]$ has a component $Q^{**}$. Recall that $Q^{**}$ is the graph obtained from $K_{1}\vee K_{6}$ and $\overline{K}_{3}\vee K_{5}$ by connecting each vertex in the part $V(\overline{K}_{3})$ to exactly two vertices in the part $V(K_{6})$. Clearly, $Q^{**}$ has an equitable partition (or regular partition, see \cite{CRS}) with 4 parts: $V(K_{1}),V(K_{6}),V(\overline{K}_{3}),V(K_{5})$. 
Then $G$ has an equitable partition with 6 parts: $V(K_{1}),V(K_{6}),V(\overline{K}_{3}),V(K_{5}),L-V(Q^{**}),R$. 
The quotient matrix of this partition is
\begin{center}$B_{7}$ =
$\left(\begin{array}{cccccc}
0&6&0&0&0&\frac{n}{2}\\
1&5&1&0&0&\frac{n}{2}\\
0&2&0&5&0&\frac{n}{2}\\
0&0&3&4&0&\frac{n}{2}\\
0&0&0&0&7&\frac{n}{2}\\
1&6&3&5&\frac{n}{2}-15&0
\end{array}\right)$.
   \end{center}
Let $h_{7}(x)$ denote the characteristic polynomial of $B_{7}$. Then $\rho(G)$ is the largest root of $h_{7}(x)=0$.
By a calculation, we have
$$h_{7}(x)=g(x)a_{7}(x)+\frac{n^{3}x}{8}-\frac{n^{4}}{16}-\frac{5xn^{2}}{4}+\frac{n^{3}}{2}
+\Theta(n^{2}),$$
where 
$$a_{7}(x)=x^{4}-9x^{3}-4x^{2}+(\frac{n}{2}+109)x
-\frac{n^{2}}{4}+2n+108.$$
When $x\geq\frac{7+\sqrt{7^{2}+n^{2}-4}}{2}>\frac{7+n}{2}$, it is easy to see that
$a_{7}(x)\geq0$ and 
$$\frac{n^{3}x}{8}-\frac{n^{4}}{16}-\frac{5xn^{2}}{4}+\frac{n^{3}}{2}
+\Theta(n^{2})>0.$$
 Note that $g(x)\geq0$ for $x\geq\frac{7+\sqrt{7^{2}+n^{2}-4}}{2}$. Hence, $h_{7}(x)>0$ for $x\geq\frac{7+\sqrt{7^{2}+n^{2}-4}}{2}$. This implies that $\rho(G)<\frac{7+\sqrt{7^{2}+n^{2}-4}}{2}$, as desired.

It remains to consider $k\geq3$ and $k\neq7$. In this case, $G[L]$ contains a component $Q^{*}_{k}$. Recall that $Q^{*}_{k}$ is the graph obtained from $K_{1}\vee K_{k-1}$ and $K_{\frac{k-1}{2}}\vee K_{\frac{k+1}{2}}$ by adding a single vertex $w$, and connecting $w$ to all vertices in the part $V(K_{\frac{k+1}{2}})$ and to $\frac{k-1}{2}$ vertices (say $u_{1},u_{2},...,u_{\frac{k-1}{2}}$) in $V(K_{k-1})$, and then connecting the remained $\frac{k-1}{2}$ vertices (say $u_{\frac{k+1}{2}},u_{\frac{k+3}{2}},...,u_{k-1}$) in $V(K_{k-1})$ to the vertices in the part $V(K_{\frac{k-1}{2}})$ by a matching of $\frac{k-1}{2}$ edges.  

Clearly, $Q^{*}_{k}$ has an equitable partition with 6 parts: $$V(K_{1}),\left\{u_{1},u_{2},...,u_{\frac{k-1}{2}}\right\},
\left\{u_{\frac{k+1}{2}},u_{\frac{k+3}{2}},...,u_{k-1}\right\},\left\{w\right\},
V(K_{\frac{k-1}{2}}),V(K_{\frac{k+1}{2}}).$$
Then $G$ has an equitable partition with 8 parts: 
$$V(K_{1}),\left\{u_{1},u_{2},...,u_{\frac{k-1}{2}}\right\},
\left\{u_{\frac{k+1}{2}},u_{\frac{k+3}{2}},...,u_{k-1}\right\},\left\{w\right\},
V(K_{\frac{k-1}{2}}),V(K_{\frac{k+1}{2}}),L-V(Q^{*}_{k}),R.$$
 The quotient matrix of this partition is
\begin{center}$B_{k}$ =
$\left(\begin{array}{cccccccc}
0&\frac{k-1}{2}&\frac{k-1}{2}&0&0&0&0&\frac{n}{2}\\
1&\frac{k-3}{2}&\frac{k-1}{2}&1&0&0&0&\frac{n}{2}\\
1&\frac{k-1}{2}&\frac{k-3}{2}&0&1&0&0&\frac{n}{2}\\
0&\frac{k-1}{2}&0&0&0&\frac{k+1}{2}&0&\frac{n}{2}\\
0&0&1&0&\frac{k-3}{2}&\frac{k+1}{2}&0&\frac{n}{2}\\
0&0&0&1&\frac{k-1}{2}&\frac{k-1}{2}&0&\frac{n}{2}\\
0&0&0&0&0&0&k&\frac{n}{2}\\
1&\frac{k-1}{2}&\frac{k-1}{2}&1&\frac{k-1}{2}&\frac{k+1}{2}&\frac{n}{2}-2k-1&0
\end{array}\right)$.
   \end{center}
Let $h_{k}(x)$ denote the characteristic polynomial of $B_{k}$. Then $\rho(G)$ is the largest root of $h_{k}(x)=0$.
By a calculation, we have
$$h_{k}(x)=g(x)a_{k}(x)+\frac{(x+k)n^{5}}{32}-\frac{n^{6}}{64}
-\frac{(k+5)xn^{4}}{16}+\frac{n^{5}}{8}+\Theta(n^{4})$$
where $$a_{k}(x)=x^{6}-(2k-5)x^{5}+(k^{2}-9k+8)x^{4}+\frac{9k^{2}+n+5-28k}{2}x^{3}
-\frac{2k^{3}+n^{2}-29k^{2}-8n+28k+23}{4}x^{2}$$
$$+\frac{n^{3}+4kn-10k^{3}+36k^{2}-10n^{2}+12n+22k-24}{8}x$$
$$-\frac{n^{4}-2kn^{3}-8n^{3}+4kn^{2}
-4k^{2}-12nk^{2}+12k^{3}+28n^{2}-68k+36n-84}{16}.$$
Clearly, $a_{k}(x)\geq n^{6}(\frac{1}{2^{6}}-o(1))>0$ for any $\frac{n}{2}\leq x\leq n$.
When $x\geq\frac{k+\sqrt{k^{2}+n^{2}-4}}{2}\geq\frac{k+n}{2}$, it is easy to see $$\frac{(x+k)n^{5}}{32}-\frac{n^{6}}{64}
-\frac{(k+5)xn^{4}}{16}+\frac{n^{5}}{8}+\Theta(n^{4})\geq\frac{k-2}{64}n^{5}+\Theta(n^{4})>0.$$
Note that $g(x)\geq0$ for $x\geq\frac{k+\sqrt{k^{2}+n^{2}-4}}{2}$. Hence, $h_{k}(x)>0$ for $x\geq\frac{k+\sqrt{k^{2}+n^{2}-4}}{2}$. This implies that $\rho(G)<\frac{k+\sqrt{k^{2}+n^{2}-4}}{2}$, as desired.
This completes the proof. \hfill$\Box$

\section{Proof of Theorem \ref{main1}}

Recall that $H_{2k+3}=K_{1}\vee P_{2k+3}$ for any $k\geq1$.
Clearly, $H_{2k+3}$ is a vertex-critical graph with $\chi(H_{2k+3})=3$.

\medskip

\f{\bf Observation 1.} $G\vee \overline{K_{2k+3}}$ is $H_{2k+3}$-free if and only if $\Delta(G)\leq k$ and $G$ is $P_{2k+3}$-free.

\medskip

\begin{lem}\label{Lemma spec-lower} Let $G\in {\rm SPEX}(n,H_{2k+3})$, where $n\geq6k$.
For odd $k\geq1$,
\begin{equation}
\rho(G)\geq\left\{
\begin{array}{lr}
\frac{k+\sqrt{k^{2}+n^{2}}}{2},&n\equiv 0(\mod4);\\
\\
\frac{k+\sqrt{k^{2}+n^{2}-4}}{2},&n\equiv 2(\mod4);\\
\\
\frac{k+\sqrt{k^{2}+n^{2}-1}}{2},&n\equiv 1,3(\mod4).
\end{array}
\right.\notag
\end{equation}
For even $k\geq2$,
\begin{equation}
\rho(G)\geq\left\{
\begin{array}{lr}
\frac{k+\sqrt{k^{2}+n^{2}}}{2},&n\equiv 0(\mod2);\\
\\
\frac{k+\sqrt{k^{2}+n^{2}-1}}{2},&n\equiv 1(\mod2).
\end{array}
\right.\notag
\end{equation}.
 \end{lem}

\f{\bf Proof:} Assume that $k\geq1$ is odd. When $n\equiv 0(\mod4)$, let $H$ be a $k$-regular graph of order $\frac{n}{2}$ such that each component of $H$ is of order at most $2k$ (such graph $H$ exists by Fact 1 in the last section). Clearly, $H$ contains no $P_{2k}$. By Observation 1,  $H\vee \overline{K_{\frac{n}{2}}}$ is $H_{2k+3}$-free. Clearly, $H\vee \overline{K_{\frac{n}{2}}}$ has an equitable partition with 2 parts: $V(H)$ and $V(\overline{K_{\frac{n}{2}}})$. Using its quotient matrix, we obtain that  $\rho(H\vee \overline{K_{\frac{n}{2}}})=\frac{k+\sqrt{k^{2}+n^{2}}}{2}$. Since $G\in {\rm SPEX}(n,H_{2k+3})$, we have
$\rho(G)\geq\rho(H\vee \overline{K_{\frac{n}{2}}})=\frac{k+\sqrt{k^{2}+n^{2}}}{2},$ as desired. For other cases, the proofs are similar. We only give the constructions of $H$. When $n\equiv 2(\mod4)$, let $H$ be a $k$-regular graph of order $\frac{n}{2}+1$ such that each component of $H$ is of order at most $2k$. When $n\equiv 1(\mod4)$, let $H$ be a $k$-regular graph of order $\frac{n-1}{2}$ such that each component of $H$ is of order at most $2k$. When $n\equiv 3(\mod4)$, let $H$ be a $k$-regular graph of order $\frac{n+1}{2}$ such that each component of $H$ is of order at most $2k$.

Assume that $k\geq2$ is even. When $n\equiv 0(\mod2)$, let $H$ be a $k$-regular graph of order $\frac{n}{2}$ such that each component of $H$ is of order at most $2k$. When $n\equiv 1(\mod2)$, let $H$ be a $k$-regular graph of order $\frac{n+1}{2}$ such that each component of $H$ is of order at most $2k$.
This completes the proof. \hfill$\Box$

\medskip

Now we are ready to prove Theorem \ref{main1}.

\medskip

\f{\bf Proof of Theorem \ref{main1}.} Let $G\in {\rm SPEX}(n,H_{2k+3})$, where $k\geq1$ and $n$ is sufficiently large. Let $\mathbf{x}=(x_{v})$ be a Perron vector of $G$ with the largest entry $1$. Let $\theta>0$ be a small constant with respect to $n$. Recall that $H_{2k+3}$ is a vertex-critical graph with $\chi(H_{2k+3})=2+1$. By Theorem \ref{skeleton},  there is a partition $V(G)=V_{1}\cup V_{2}$ such that $||V_{i}|-\frac{n}{2}|<\theta n$ for any  $1\leq i\leq 2$, and $d_{V_{i}}(v)< 2k+3$ and $d_{G}(v)\geq n-|V_{i}|$ for any $v\in V_{i}$. Moreover, $1-\theta<x_{v}\leq1$ for any $v\in V(G)$. 
Very similar to Claim 5 of Theorem \ref{skeleton}, the following Claim 1 can be proved.

\medskip

\f{\bf Claim 1.}
Let $1\leq \ell\leq 2$ be fixed. For $1\leq i\neq\ell\leq2$, assume that  $u_{1},u_{2},...,u_{2k+3}\in V_{i}$. Then there are $2k+3$ vertices in $V_{\ell}$ which are adjacent to all the vertices $u_{1},u_{2},...,u_{2k+3}$ in $G$.

\medskip

\medskip

\f{\bf Claim 2.}
$\Delta(G[V_{i}])\leq k$ for any $1\leq i\leq 2$. (This implies that any $v\in V_{i}$ has at most $k$ non-neighbors in $V_{3-i}$ as $d_{G}(v)\geq n-|V_{i}|$.)

\medskip

\f{\bf Proof of Claim 2.} Suppose not. Without loss of generality, let $v\in V_{1}$ with $d_{V_{1}}(v)\geq k+1$. Then $v$ has $k+1$ neighbors in $V_{1}$, say $u_{1},u_{2},...,u_{k+1}$. By Claim 1, $v,u_{1},u_{2},...,u_{k+1}$ have at least $k+2$ common neighbors in $V_{2}$, say $u_{k+2},u_{k+3},...,u_{2k+3}$. Clearly, the subgraph of $G$ induced by $\left\{v,u_{1},u_{2},...,u_{2k+3}\right\}$ contains a copy of $H_{2k+3}$, a contradiction. Hence, $\Delta(G[V_{i}])\leq k$ for any $1\leq i\leq 2$.
 This finishes the proof of Claim 2. \hfill$\Box$

\medskip

\f{\bf Claim 3.} For a constant $C>0$ (with respect to $n$),
let $G'$ be the graph obtained from $G$ by deleting at most $C$ edges. Then $\Delta(G'[V_{1}])+\Delta(G'[V_{2}])\geq k$.

\medskip

\f{\bf Proof of Claim 3.} Let $d_{i}=\Delta(G'[V_{i}])$ for $1\leq i\leq2$. By Lemma \ref{subgraph} and Lemma \ref{loop}, we have $$\rho(G')\leq\frac{d_{1}+d_{2}+\sqrt{(d_{1}-d_{2})^{2}+4|V_{1}||V_{2}|}}{2}.$$
Let $\mathbf{x}^{T}$ denote the transpose of $\mathbf{x}$. Since $1-\theta<x_{v}\leq1$ for any $v\in V(G)$, we have
\begin{equation}
\begin{aligned}
\rho(G)&=\frac{\mathbf{x}^{T}A(G)\mathbf{x}}{\mathbf{x}^{T}\mathbf{x}}\\
&\leq\frac{\mathbf{x}^{T}A(G')\mathbf{x}+2C}{\mathbf{x}^{T}\mathbf{x}}\\
&\leq\rho(G')+\frac{2C}{n(1-\theta)^{2}}\\
&\leq\frac{d_{1}+d_{2}+\sqrt{(d_{1}-d_{2})^{2}+4|V_{1}||V_{2}|}}{2}+\frac{2C}{n(1-\theta)^{2}}\\
&\leq\frac{d_{1}+d_{2}+\sqrt{(d_{1}-d_{2})^{2}+n^{2}}}{2}+\frac{8C}{n}.
\end{aligned}\notag
\end{equation}
By Lemma \ref{Lemma spec-lower}, we have $\rho(G)\geq\frac{k+\sqrt{k^{2}+n^{2}-4}}{2}$. Thus,
$$\frac{d_{1}+d_{2}+\sqrt{(d_{1}-d_{2})^{2}+n^{2}}}{2}+\frac{8C}{n}\geq\frac{k+\sqrt{k^{2}+n^{2}-4}}{2},$$
implying that $d_{1}+d_{2}\geq k$ for sufficiently large $n$.
 \hfill$\Box$

\medskip

Without loss of generality, assume that $\Delta(G[V_{1}])\geq\Delta(G[V_{2}])$.

\medskip

\f{\bf Claim 4.} $\Delta(G[V_{1}])=k$.

\medskip

\f{\bf Proof of Claim 4.} Since $\Delta(G[V_{1}])+\Delta(G[V_{2}])\geq k$ by Claim 3 and $\Delta(G[V_{1}])\geq\Delta(G[V_{2}])$, we have $\Delta(G[V_{1}])\geq\lceil\frac{k}{2}\rceil$. If $k=1$, then $\Delta(G[V_{1}])=1$, as desired. Now assume that $k\geq2$ in the following.

Now we prove $\Delta(G[V_{1}])\geq k-1$. Suppose not. Then $\Delta(G[V_{1}])\leq k-2$. Since $\Delta(G[V_{1}])\geq\lceil\frac{k}{2}\rceil$, we have $k\geq4$. Let $u$ be a vertex in $V_{1}$ such that $d_{V_{1}}(u)=\Delta(G[V_{1}])$. Then $u$ has at least $\lceil\frac{k}{2}\rceil$ neighbors in $V_{1}$, say $u_{1},u_{2},...,u_{\lceil\frac{k}{2}\rceil}$. By Claim 2, $u,u_{1},u_{2},...,u_{\lceil\frac{k}{2}\rceil}$ are adjacent to all the vertices in $T\subseteq V_{2}$, where $|T|\geq|V_{2}|-k\lceil\frac{k+2}{2}\rceil$. If $G[T]$ has $\lceil\frac{k+2}{2}\rceil(1+k+k(k-1))$ vertices of degree at least 2, then $G[T]$ has $\lceil\frac{k+2}{2}\rceil$ vertex-disjoint copies of $P_{3}$. Clearly, the subgraph induced by $u_{1},u_{2},...,u_{\lceil\frac{k}{2}\rceil}$ and the vertices in these $\lceil\frac{k+2}{2}\rceil$ vertex-disjoint copies of $P_{3}$, contains a  path of order $\lceil\frac{k}{2}\rceil+3\lceil\frac{k+2}{2}\rceil\geq2k+3$. Adding the vertex $u$, there is a copy of $H_{2k+3}$ in $G$, a contradiction. Thus, $G[T]$ has less than $\lceil\frac{k+2}{2}\rceil(1+k+k(k-1))=\lceil\frac{k+2}{2}\rceil(k^{2}+1)$ vertices of degree at least 2. Let $G_{1}$ be the graph obtained from $G$ by deleing the edges incident with vertices of degree $\geq2$ in $V_{2}$. Then $\Delta(G_{1}[V_{1}])\leq k-2$ and $\Delta(G_{1}[V_{2}])\leq 1$. However, we have deleted at most $k^{2}\lceil\frac{k+2}{2}\rceil+k\lceil\frac{k+2}{2}\rceil(k^{2}+1)$ edges from $G$. By Claim 3, we have $\Delta(G_{1}[V_{1}])+\Delta(G_{1}[V_{2}])\geq k$, a contradiction. Hence, we must have $\Delta(G[V_{1}])\geq k-1$.

Now we prove $\Delta(G[V_{1}])=k$. Suppose not. Then $\Delta(G[V_{1}])=k-1$. We will obtain a contradiction by three cases as follows.

\medskip

\f{\bf Case 1.} $k=2$.

\medskip

In this case, $\Delta(G[V_{2}])\leq\Delta(G[V_{1}])=1$.
 By Lemma \ref{subgraph} and Lemma \ref{loop}, we have $$\rho(G)\leq\frac{1+1+\sqrt{(1-1)^{2}+4|V_{1}||V_{2}|}}{2}\leq\frac{2+\sqrt{n^{2}}}{2}.$$
But $$\rho(G)\geq\frac{2+\sqrt{2^{2}+n^{2}-1}}{2}=\frac{2+\sqrt{n^{2}+3}}{2}$$
 by Lemma \ref{Lemma spec-lower} as $k=2$, a contradiction.

\medskip

\f{\bf Case 2.} $k=3$.

\medskip

In this case, $\Delta(G[V_{2}])\leq\Delta(G[V_{1}])=2$. If $\Delta(G[V_{2}])\leq1$,
by Lemma \ref{subgraph} and Lemma \ref{loop}, we have $$\rho(G)\leq\frac{2+1+\sqrt{(2-1)^{2}+4|V_{1}||V_{2}|}}{2}\leq\frac{3+\sqrt{1+n^{2}}}{2}.$$
But $$\rho(G)\geq\frac{3+\sqrt{3^{2}+n^{2}-4}}{2}=\frac{3+\sqrt{n^{2}+5}}{2}$$ by Lemma \ref{Lemma spec-lower} as $k=3$, a contradiction.

It remains that $\Delta(G[V_{2}])=\Delta(G[V_{1}])=2$. Let $u$ be a vertex in $V_{1}$ such that $d_{V_{1}}(u)=2$. Then $u$ has two neighbors in $V_{1}$, say $u_{1},u_{2}$. By Claim 2, $u,u_{1},u_{2}$ are adjacent to all the vertices in $T\subseteq V_{2}$, where $|T|\geq|V_{2}|-9$. If $G[T]$ has $3(1+3+3(3-1))$ vertices of degree at least 2, then $G[T]$ has $3$ vertex-disjoint copies of $P_{3}$. Clearly, the subgraph induced by $u_{1},u_{2}$ and the vertices in these $3$ vertex-disjoint copies of $P_{3}$, contains a  path of order $11\geq9$. Adding the vertex $u$, there is a copy of $H_{9}$ in $G$, a contradiction. Thus, $G[T]$ has less than $3(1+3+3(3-1))=30$ vertices of degree at least 2. Thus, $G[V_{2}]$ will have maximum degree 1 after deleting at most $9\cdot3+30\cdot3=117$ edges. Recall that $\Delta(G[V_{2}])=2$. Similarly, $G[V_{1}]$ will also have maximum degree 1 after deleting at most $117$ edges. Let $G_{2}$ be the graph obtained from $G$ by deleing the $\leq117\cdot2=234$ edges inside $V_{1}$ and $V_{2}$. Then $\Delta(G_{2}[V_{1}])=\Delta(G_{2}[V_{2}])=1$. By Claim 3, we have $1+1\geq 3$, a contradiction.

\medskip

\f{\bf Case 3.} $k\geq4$.

\medskip

In this case, $\Delta(G[V_{2}])\leq\Delta(G[V_{1}])=k-1$. Let $u$ be a vertex in $V_{1}$ such that $d_{V_{1}}(u)=\Delta(G[V_{1}])$. Then $u$ has $k-1$ neighbors in $V_{1}$, say $u_{1},u_{2},...,u_{k-1}$. By Claim 2, $u,u_{1},u_{2},...,u_{k-1}$ are adjacent to all the vertices in $T\subseteq V_{2}$, where $|T|\geq|V_{2}|-k^{2}$. If $G[T]$ has $4$ vertex-disjoint edges, say $v_{1}w_{1},v_{2}w_{2},v_{3}w_{3},v_{4}w_{4},$ select other $k-4$ vertices in $T$, say $v_{5},v_{6},...,v_{k}$. Clearly, the subgraph induced by $u_{1},u_{2},...,u_{k-1}$ and the vertices $v_{1},w_{1},v_{2},w_{2},v_{3},w_{3},v_{4},w_{4},v_{5},v_{6},...,v_{k}$, contains a  path of order $2k+3$. Adding the vertex $u$, there is a copy of $H_{2k+3}$ in $G$, a contradiction. Thus, $G[T]$ has at most $3$ vertex-disjoint edges. Since $\Delta(G[T])\leq k$, we have $e(G[T])\leq6k$. It follows that $e(G[V_{2}])\leq k^{3}+6k$. Let $G_{3}$ be the graph obtained from $G$ by deleing all the edges inside $V_{2}$. Then $\Delta(G_{3}[V_{1}])= k-1$ and $\Delta(G_{3}[V_{2}])=0$. By Claim 3, we have $k-1=\Delta(G_{3}[V_{1}])+\Delta(G_{3}[V_{2}])\geq k$, a contradiction.
 This finishes the proof of Claim 4. \hfill$\Box$

\medskip

\f{\bf Claim 5.} $G[V_{1}]$ contains no $P_{2k+3}$. Thus, each component of $G[V_{1}]$ has order at most $f(k)$. Recall that $f(k)=\sum_{0\leq i\leq 2k+1}k^{i}$.

\medskip

\f{\bf Proof of Claim 5.} If $G[V_{1}]$ contains a path $v_{1}v_{2}\cdots v_{2k+3}$, by Claim 1, the vertices in this path have a common neighbor in $V_{2}$. Thus, $H_{2k+3}$ arises in $G$, a contradiction. Hence $G[V_{1}]$ contains no $P_{2k+3}$. 
This finishes the proof of Claim 5. \hfill$\Box$

\medskip

\f{\bf Claim 6.} $e(G[V_{2}])=0$. Each vertex in $V_{1}$ is adjacent to all the vertices in $V_{2}$.

\medskip

\f{\bf Proof of Claim 6.} We first show that $e(G[V_{2}])\leq k^{3}+k^{2}+2k$. By Claim 4, there is a vertex $u$ in $V_{1}$, which has $k$ neighbors in $V_{1}$, say $u_{1},u_{2},...,u_{k}$. By Claim 2, $u,u_{1},u_{2},...,u_{k}$ are adjacent to all the vertices in $T\subseteq V_{2}$, where $|T|\geq|V_{2}|-k(k+1)$.
If $G[T]$ has $2$ vertex-disjoint edges, say $v_{1}w_{1},v_{2}w_{2}$, select other $k-1$ vertices in $T$, say $v_{3},v_{4},...,v_{k+1}$. Clearly, the subgraph induced by $u_{1},u_{2},...,u_{k}$ and the vertices $v_{1},w_{1},v_{2},w_{2},v_{3},v_{4},...,v_{k+1}$, contains a  path of order $2k+3$. Adding the vertex $u$, there is a copy of $H_{2k+3}$ in $G$, a contradiction. Thus, $G[T]$ has at most $1$ vertex-disjoint edge. Since $\Delta(G[T])\leq k$, we have $e(G[T])\leq2k$. It follows that $e(G[V_{2}])\leq k^{2}(k+1)+2k=k^{3}+k^{2}+2k$.

Now we show $e(G[V_{2}])=0$. Suppose not. Let $v_{0}w_{0}$ be an edge inside $V_{2}$. By Claim 2, $v_{0},w_{0}$ are adjacent to all the vertices in $S\subseteq V_{1}$, where $|S|\geq|V_{1}|-2k$. Since $n$ is large, we can require $\frac{|V_{1}|}{f(k)}\geq2k+1+2(k^{3}+k^{2}+2k+1)$. This implies that there are at least $1+2(k^{3}+k^{2}+2k+1)$ components of $G[V_{1}]$, say $Q_{1},Q_{2},...,Q_{1+2(k^{3}+k^{2}+2k+1)}$, such that $V(Q_{i})\subseteq S$ for any $1\leq i\leq1+2(k^{3}+k^{2}+2k+1)$. If for some $1\leq i<j\leq1+2(k^{3}+k^{2}+2k+1)$, both  $Q_{i}$ and $Q_{j}$ contain $P_{k+1}$, then the subgraph induced by $w_{0}$ and the vertices in $V(Q_{i})\cup V(Q_{j})$, contains a  path of order $2k+3$. Adding the vertex $v_{0}$, there is a copy of $H_{2k+3}$ in $G$, a contradiction. Thus, there is at most one component, say $Q_{1+2(k^{3}+k^{2}+2k+1)}$, which contains a path of order $k+1$. That is, $Q_{i}$ contains no $P_{k+1}$ for any $1\leq i\leq 2(k^{3}+k^{2}+2k+1)$, which implies that $Q_{i}$ is not $k$-regular. Hence, there is a vertex $z_{i}$ with degree less than $k$ in $Q_{i}$  for any $1\leq i\leq 2(k^{3}+k^{2}+2k+1)$. Let $G_{4}$ be the graph obtained from $G$ by deleing all the edges inside $V_{2}$, and adding the edges $z_{2i-1}z_{2i}$ for any $1\leq i\leq k^{3}+k^{2}+2k+1$. Since $Q_{i}$ contains no $P_{k+1}$ for any $1\leq i\leq 2(k^{3}+k^{2}+2k+1)$, $G_{4}[V_{1}]$ contains no $P_{2k+3}$. Note that $\Delta(G_{4}[V_{1}])=k$. Then $G_{4}$ is $H_{2k+3}$-free by Observation 1. Since $1-\theta\leq x_{w}\leq1$ for any $w\in V(G)$, it is easy to see
\begin{equation}
\begin{aligned}
&\mathbf{x}^{T}(\rho(G_{4})-\rho(G))\mathbf{x}\\
&\geq\mathbf{x}^{T}(A(G_{4})-A(G))\mathbf{x}\\
&=2(\sum_{uv\in E(G_{4})-E(G)}x_{u}x_{v})-2(\sum_{uv\in E(G)-E(G_{4})}x_{u}x_{v})\\
&\geq2(1-\theta)^{2}(k^{3}+k^{2}+2k+1)-2(k^{3}+k^{2}+2k)\\
&>0~(requiring~\theta<1-\sqrt{\frac{k^{3}+k^{2}+2k}{k^{3}+k^{2}+2k+1}}).
\end{aligned}\notag
\end{equation}
It follows that $\rho(G_{4})>\rho(G)$. But this contradicts that $G\in {\rm SPEX}(n,H_{2k+3})$. Hence, $e(G[V_{2}])=0$.

Since $G[V_{1}]$ contains no $P_{2k+3}$ by Claim 5, by Lemma \ref{subgraph} and Observation 1, we must have that each vertex in $V_{1}$ is adjacent to all the vertices in $V_{2}$. This finishes the proof of Claim 6. \hfill$\Box$

\medskip

\f{\bf Claim 7.} $G[V_{1}]$ is $k$-regular or nearly $k$-regular.

\medskip

\f{\bf Proof of Claim 7.} Suppose not. We can choose the union of some components of $G[V_{1}]$, say $Q$, such that $4k+3\leq|Q|\leq2f(k)$ and $e(Q)\leq \lfloor\frac{k}{2}|Q|\rfloor-1$. Let $G_{5}$ be a graph obtained from $G$ by deleing all the edges of $Q$, and embedding a $k$-regular or nearly $k$-regular $P_{2k+3}$-free graph in $V(Q)$. Clearly, $G_{5}$ is $H_{2k+3}$-free by Observation 1. Since $1-\theta\leq x_{w}\leq1$ for any $w\in V(G)$, we have
\begin{equation}
\begin{aligned}
&\mathbf{x}^{T}(\rho(G_{5})-\rho(G))\mathbf{x}\\
&\geq\mathbf{x}^{T}(A(G_{5})-A(G))\mathbf{x}\\
&=2(\sum_{uv\in E(G_{5})-E(G)}x_{u}x_{v})-2(\sum_{uv\in E(G)-E(G_{5})}x_{u}x_{v})\\
&\geq2(1-\theta)^{2}\lfloor\frac{k}{2}|Q|\rfloor-2(\lfloor\frac{k}{2}|Q|\rfloor-1)\\
&\geq2((1-\theta)^{2}kf(k)-kf(k)+1)~(|Q|\leq2f(k))\\
&>0~(requiring~\theta<1-\sqrt{\frac{kf(k)-1}{kf(k)}}).
\end{aligned}\notag
\end{equation}
It follows that $\rho(G_{5})>\rho(G)$. But this contradicts that $G\in {\rm SPEX}(n,H_{2k+3})$. Hence, $G[V_{1}]$ is $k$-regular or nearly $k$-regular.
 This finishes the proof of Claim 7. \hfill$\Box$

\medskip

Now we prove the theorem by cases.

\medskip

\f{\bf Case 1.}  $k\geq2$ is even.

\medskip

Since $k$ is even, we see that $G[V_{1}]$ is $k$-regular. Then $G$ has an equitable partition: $V_{1}$ and $V_{2}$. Using quotient matrix, we have $\rho(G)=\frac{k+\sqrt{k^{2}+4|V_{1}||V_{2}|}}{2}$. By Lemma \ref{Lemma spec-lower}, we have $\rho(G)\geq\frac{k+\sqrt{k^{2}+n^{2}-1}}{2}$. Thus,
$$\frac{k+\sqrt{k^{2}+n^{2}-1}}{2}\leq\frac{k+\sqrt{k^{2}+4|V_{1}||V_{2}|}}{2}.$$
It follows that $|V_{1}|=\lfloor\frac{n}{2}\rfloor$ or $|V_{1}|=\lceil\frac{n}{2}\rceil$.

\medskip

\f{\bf Case 2.}  $k\geq1$ is odd.

\medskip
 
 In this case, we will prove
\begin{equation}
|V_{1}|=\left\{
\begin{array}{lr}
\frac{n}{2},&n\equiv 0(\mod4);\\
\\
\frac{n-1}{2},&n\equiv 1(\mod4);\\
\\
\frac{n}{2},&k=1~and~n\equiv 2(\mod4);\\
\\
\frac{n}{2}-1~or~\frac{n}{2}+1,&k\geq3~and~n\equiv 2(\mod4);\\
\\
\frac{n+1}{2},&n\equiv 3(\mod4).\\
\end{array}
\right.\notag
\end{equation}
 Let
\begin{center}$B$=
$\left(\begin{array}{cc}
k&|V_{2}|\\
|V_{1}|&0
\end{array}\right)$.
   \end{center}
   By Lemma \ref{loop}, we have $\rho(G)\leq\rho(B)=\frac{k+\sqrt{k^{2}+4|V_{1}||V_{2}|}}{2}$. Note that the equality can not hold if $G[V_{1}]$ is not $k$-regular.
   By Lemma \ref{Lemma spec-lower},
   \begin{equation}
\rho(G)\geq\left\{
\begin{array}{lr}
\frac{k+\sqrt{k^{2}+n^{2}}}{2},&n\equiv 0(\mod4);\\
\\
\frac{k+\sqrt{k^{2}+n^{2}-4}}{2},&n\equiv 2(\mod4);\\
\\
\frac{k+\sqrt{k^{2}+n^{2}-1}}{2},&n\equiv 1,3(\mod4).
\end{array}
\right.\notag
\end{equation}

\medskip

\f{\bf Subcase 2.1.} $n\equiv 0,1,3(\mod4)$.

\medskip

In this case, it is easy to check that $|V_{1}|=\lfloor\frac{n}{2}\rfloor$ or $\lceil\frac{n}{2}\rceil$, and $G[V_{1}]$ must be $k$-regular. Then $|L|$ must be displayed in the theorem.

\medskip

\f{\bf Subcase 2.2.} $n\equiv 2(\mod4)$.

\medskip

Using Lemma \ref{Lemma spec-lower}, we have
$$\frac{k+\sqrt{k^{2}+n^{2}-4}}{2}\leq\rho(G)\leq\frac{k+\sqrt{k^{2}+4|V_{1}||V_{2}|}}{2}.$$
It follows that 
$$\frac{n}{2}-1\leq|V_{1}|\leq\frac{n}{2}+1.$$ 
If $|V_{1}|=\frac{n}{2}-1$ or $\frac{n}{2}+1$,  then $G[V_{1}]$ is $k$-regular. Let $g(x)=x^{2}-kx-(\frac{n^{2}}{4}-1)$. Then $\rho(G)=\frac{k+\sqrt{k^{2}+n^{2}-4}}{2}$ is the largest root of $g(x)=0$.
 
 \medskip
 
 \f{\bf Subcase 2.2.1.} $k=1$.
  
  \medskip
  
  In this case, we will show that $|V_{1}|=\frac{n}{2}$. In fact, If $|V_{1}|=\frac{n}{2}-1$ or $\frac{n}{2}+1$,  then as above, $\rho(G)=\frac{1+\sqrt{n^{2}-3}}{2}$ is the largest root of $g(x)=0$, where $g(x)=x^{2}-x-(\frac{n^{2}}{4}-1)$. If $|V_{1}|=\frac{n}{2}$, then $G=(K_{1}\cup M_{\frac{n}{2}-1})\vee\overline{K_{\frac{n}{2}}}$, where $M_{\frac{n}{2}-1}$ denotes a matching of order $\frac{n}{2}-1$. Clearly, $G$ has an equitable partition with 3 parts: $V(K_{1}),V(M_{\frac{n}{2}-1}),V(\overline{K_{\frac{n}{2}}})$. The quotient matrix is
\begin{center}$B_{1}$ =
$\left(\begin{array}{ccc}
0&0&\frac{n}{2}\\
0&1&\frac{n}{2}\\
1&\frac{n-2}{2}&0\\
\end{array}\right)$.
   \end{center}
Let $h_{1}(x)$ denote the characteristic polynomial of $B_{1}$. Then $\rho(G)$ is the largest root of $h_{1}(x)=0$. By a simple calculation, we have $h_{1}(x)=x^{3}-x^{2}-\frac{n^{2}}{4}x+\frac{n}{2}$,
 and
$$h_{1}(x)=xg(x)+\frac{n}{2}-x.$$
 Clearly, $h_{1}(\frac{1+\sqrt{n^{2}-3}}{2})<0$, since $g(\frac{1+\sqrt{n^{2}-3}}{2})=0$ and $\frac{1+\sqrt{n^{2}-3}}{2}>\frac{n}{2}$. This implies that $\rho(G)>\frac{1+\sqrt{n^{2}-3}}{2}$. Consequently, we must have $|V_{1}|=\frac{n}{2}$.

\medskip
 
 \f{\bf Subcase 2.2.2.} $k\geq3$.
  
  \medskip

In this case, we will prove $|V_{1}|=\frac{n}{2}-1$ or $\frac{n}{2}+1$. Suppose not. Then $|V_{1}|=\frac{n}{2}$ by above.  Clearly, $G[V_{1}]$ is nearly $k$-regular. Recall that $G[V_{1}]$ is $P_{2k+3}$-free. Hence, $G$ is obtained from the Tur\'{a}n graph $T(n,2)$ with parts $V_{1}$ and $V_{2}$, by embedding a graph from $\mathcal{G}_{\frac{n}{2},k}$ into $V_{1}$. But then $\rho(G)<\frac{k+\sqrt{k^{2}+n^{2}-4}}{2}$ by Theorem \ref{spec compare}. This contradicts the fact that $\rho(G)\geq\frac{k+\sqrt{k^{2}+n^{2}-4}}{2}$ by Lemma \ref{Lemma spec-lower}. Hence, we must have $|V_{1}|=\frac{n}{2}-1$ or $\frac{n}{2}+1$.
This completes the proof. \hfill$\Box$

\medskip

\medskip

\f{\bf Declaration of competing interest}

\medskip

There is no conflict of interest.

\medskip

\f{\bf Data availability statement}

\medskip

No data was used for the research described in the article.

\end{document}